\documentclass[sn-mathphys,Numbered]{sn-jnl}

\usepackage{graphicx}%
\usepackage{multirow}%
\usepackage{amsmath,amssymb,amsfonts}%
\usepackage{amsthm}%
\usepackage{mathrsfs}%
\usepackage[title]{appendix}%
\usepackage{xcolor}%
\usepackage{textcomp}%
\usepackage{manyfoot}%
\usepackage{booktabs}%
\usepackage{algorithm}%
\usepackage{algorithmicx}%
\usepackage{algpseudocode}%
\usepackage{listings}%

\theoremstyle{thmstyleone}%
\newtheorem{theorem}{Theorem}
\theoremstyle{thmstyletwo}%
\newtheorem{lemma}{Lemma}
\newtheorem{corollary}{Corollary}
\theoremstyle{thmstylethree}%
\def\sgn{{sgn\,}}
\begin{document}
	
	\title[Reconstruction of a function from its spherical transforms ]{Reconstruction of a function from its two spherical Radon transforms with the centers on a plane}
	
	
	\author*{\fnm{Rafik} \sur{Aramyan}}\email{rafikaramyan@yahoo.com}

	\affil{\orgdiv{Institute of Mathematics of National Academy of Sciences}, \orgaddress{{Bagramian 24/5},{Yerevan}, \postcode{0019}, \country{Armenia}}}

	\footnotetext{The work was supported by the Science Committee of RA, in the frames of the research project 21AG-1A045.}
	
	\abstract{Hyperplane is a set of non-injectivity of the spherical Radon transform (SRT) in the space of 
		continuous functions in $\mathbf R^{d}$. In this article, for the reconstruction of
		$f\in C(\mathbf R^3)$ (the support can be non-compact), using the spherical Radon transform over spheres centered on a plane, the injectivity of the so-called two data spherical Radon transform is considered.  An inversion formula of the transform that uses the local data of the spherical integrals to reconstruct the unknown function is presented. Such inversions are the mathematical base of modern modalities of imaging, such as Thermo and photoacoustic tomography and radar imaging, and have theoretical significance in many areas of mathematics.}

	
	
	
	\maketitle
	
	\section{Introduction}\label{sec1}

	\noindent For a continuous function $f \in
C(\mathbf R^d)$, the spherical mean $Mf(P, t)$  is defined by the formula

\begin{equation}\label{1.0}
	Mf(P,t)=\frac{1}{|\mathbf S^{d-1}|}\int_{\mathbf S^{d-1}}f(P+t\omega)\,d\omega,
\end{equation}
where ${\mathbf S^{d-1}}$ is the unit sphere in $\mathbf R^{d}$ (centered at the origin) and we integrate with respect to the surface
Lebesgue measure on ${\mathbf S^{d-1}}$ normalized by $|\mathbf S^{d-1}|$ (the area of the unit sphere). $Mf(P, t)$ is the integral (mean)
of f over the sphere $S(P, t)$ of radius $t > 0$  centered at $P$. The spherical Radon transform (SRT) takes function $f \in
C(\mathbf R^d)$ to its spherical means $Mf(P, t)$ on $L\times[0,\infty]$ ($l\subset{\mathbf S^{d-1}}$ is a hypersurface, set of detectors). The main problems are finding which $L$ is a set of injectivity (the data collected by detectors along $L$ sufficient for unique reconstruction of $f$) and finding the inversion formulas.

\noindent The inversion of SRT is required in mathematical models of modern modalities of imaging, such as Thermoacoustic (TAT), Photoacoustic tomography  and radar imaging, and has theoretical significance. (\cite{Amm2}, \cite{KuKu}, \cite{Kru}, \cite{Nat}, \cite{Pal}).

\noindent Inversion formulas of SRT are known for different sets (which are sets of injectivity) of detectors (\cite{AgQu}-\cite{AmKu}, \cite{An},\cite{Ara19}, \cite{Be}-\cite{FiHa}, \cite{Hal1}-\cite{No}, \cite{Pal}-\cite{XW}).

\noindent Hyperplane is a set of non-injectivity of SRT in the space of  continuous functions in $\mathbf R^{d}$ (see also \cite{AgQu}, \cite{CH}).
As usual, the injectivity of SRT in the space of compactly supported functions is considered (\cite{AgQu}-\cite{AmKu}, \cite{FiRa},\cite{FiHa}).

\noindent In this article, an additional condition for the reconstruction of
an unknown function $f\in C(\mathbf R^3)$ (the support can be non-compact) using SRT over spheres centered on a plane is found (Problem 1).

\noindent It is proved that this problem is equivalent to the injectivity of a so-called two data
spherical Radon transform.  Also, we present an inversion formula of the transform which uses the local data of the spherical integrals to reconstruct the unknown function. Here, we apply the consistency method which was already used in 2D (\cite{Ara19}) and 3D (\cite{Ara231}) to invert SRT in the class  of smooth functions supported on one side of a line and on one side of a plane respectively (see also \cite{Ara10}, \cite{Ara21}). A similar problem was already solved in 2D (\cite{Ara232}). Several other important  generalized Radon transforms in imaging, including weighted Radon transforms, were considered in  (\cite{Gin}, \cite{Nov}, \cite{P}). In (\cite{Zan}) a combination of the acoustic pressure data and its normal derivative is considered.

\noindent We consider the following transform  over spheres centered on the $\{z=0\}$-plane. Below, the point $P=(p,q,0)$ is identified with $(p,q)$. Given a continuous function 
$f\in\textit{C}(\mathbf R^3)$ for $(P,t)\in \{z=0\}\times(0,\infty)$ we define 
\begin{equation}\label{1.2}
	Tf(P,t)=Tf(p,q,t)=	\frac{1}{4\pi}\int_{\mathbf S^{2}}(1+i \cos\theta)f((p,q,0)+t\omega)\,d\omega, 
\end{equation}
here $(\theta,\varphi)=\omega\in\mathbf S^{2}$ is the spherical coordinates of $\omega$ ($\theta\in[0,\pi]$  is the polar angle measured from the $z$-axis) and we integrate with
respect to the surface Lebesgue measure on ${\mathbf S^{2}}$.

\noindent We denote by ${C}^{1}(\mathbf R^d)$ the space of all functions with continuous partial derivatives of first order and by  ${ C}^{\infty}(\mathbf R^d)$ the space of all functions with continuous partial derivatives of all orders.

\noindent \emph{The main results that we will prove in the following sections}.

\begin{theorem}\label{1} Let $D\subset\{z=0\}$ be a planar domain (an open subset). The transform (two data
	spherical Radon transform)
	
	\begin{equation}\label{1.3}
		f \to 	Tf(p,q,t)\,\,\,\,\emph{for}\,\,\,\,(p,q,t)\in D\times(0,\infty)\end{equation}
	is invertible in ${C}^{1}(\mathbf R^{3})$ .
\end{theorem}

\noindent For the polar coordinates $(t,\omega)=(t,\theta,\varphi)$ of $(x,y,z)\in\mathbf R^{3}$ with respect to a fixed $P=(p,q,0)$ we have 
\begin{equation}\label{1.3.1}\begin{cases}
		x=p+t \sin\theta\cos\varphi\\
		y=q+t\sin\theta\sin\varphi \\
		z=t \cos\theta.\end{cases}
\end{equation}
here $(\theta,\varphi)$ is the spherical coordinates of $\omega\in\mathrm{S}^{2}$.

\noindent For any $(p,q,t)$	by $f_{(p,q,t)}(\theta,\varphi)$ we denote the restriction of $f$ onto the sphere $S(p,q,t)$ with center $P=(p,q,0)$ and radius $t > 0$. 

\noindent For $f\in{C}^{\infty}(\mathbf R^3)$	the spherical harmonic expansion of $f_{(p,q,t)}(\theta,\varphi)$ we write as follows: 
\begin{eqnarray}\label{1.4}
	f_{(p,q,t)}(\theta,\varphi)=Mf(p,q,t)+\sum_{n=1}^{\infty} \left[a_{0n}(p,q,t) \,P_n(\cos\theta)+\right.\\\left.\sum_{m=1}^n (a_{mn}(p,q,t)\cos(m\varphi) + b_{mn}(p,q,t)\sin(m\varphi))P_{nm}(\cos\theta)\right],\nonumber\end{eqnarray}
where $Mf(p,q,t)=a_{00}(p,q,t)$ is the the spherical mean  of $f$ over $S(p,q,t)$ and ($n\geq0$ and $1\leq m\leq n$)
\begin{equation}\label{1.4.1}\begin{cases}
		a_{0n}(p,q,t)=\frac{2n+1}{4\pi}\int_{\mathbf S^{2}}f_{(p,q,t)}(\theta,\varphi)P_n(\cos\theta)\,d\omega,\\
		a_{mn}(p,q,t)=\frac{2n+1}{2\pi}\frac{(n-m)!}{(n+m)!}\int_{\mathbf S^{2}}f_{(p,q,t)}(\theta,\varphi)P_{n,m}(\cos\theta)\,\cos(m\varphi)\,d\omega,
		\\
		b_{mn}(p,q,t)=\frac{2n+1}{2\pi}\frac{(n-m)!}{(n+m)!}\int_{\mathbf S^{2}}f_{(p,q,t)}(\theta,\varphi)P_{n,m}(\cos\theta)\,\sin(m\varphi)\,d\omega.\end{cases}\end{equation}
Here $P_n$ are the Legendre polynomials and
$P_{nm}=(-1)^m P_n^m$ ($n\geq0$ and $1\leq m\leq n$),
where $P_n^m$ are the associated Legendre polynomials. Later on, we will use the short form of the integrals, for example
$$\int_{\mathbf S^{2}}f_{(p,q,t)}(\theta,\varphi)P_n(\cos\theta)\,d\omega=\int_{S(p,q,t)}f\cdot P_n(\cos\theta)\,d\omega.$$

\noindent From (\ref{1.2}) it follows that ($P_1(\cos\theta)=\cos\theta$)
\begin{equation}Re(Tf(p,q,t)) =Mf(p,q,t),\,\,\,\,\,\,\,\,Im(Tf(p,t)) =\frac{1}{3}a_{01}(p,q,t)\label{1.5}\end{equation}
where $Re(Tf(p,q,t))$ is the real and $Im(Tf(p,t))$ is the imaginary parts of $Tf(p,q,t)$, $Mf(p,q,t)$ and $a_{01}(p,q,t)$ are the Fourier first coefficients of the restriction of $f$ onto $S(p,q,t)$.

\noindent Thus, from Theorem 1 it follows that: for reconstructing an unknown function $f$ (the support can be non-compact) it is sufficient to have the Fourier first two coefficients ($Mf(p,q,t)$ and $a_{01}(p,q,t)$) of the restrictions of $f$ onto $S(p,q,t)$, $(p,q, t)\in D\times[0,\infty)$.

\noindent The primary goal of this article is to present an inversion formula for (\ref{1.3}). 

\noindent We construct two sequences of standard (they do not depend on $f$)
polynomials $C_{2n,2i}(t)$ and $C_{2n+1,2i}(t)$ ($n,i$ are integers and $0\leq i\leq n$) defined on $[0,1]$ of degree $2n+2i$ and $2n+2i+1$ respectively (see (\ref{7.2}) below).

\noindent Given a smooth function $f\in{\textit{C}}^{\infty}(\mathbf R^3)$. Let $Mf$ be the spherical mean and $a_{01}$ be the Fourier first coefficient of the restriction of $f$ onto sphere $S(P,t)$ with center $P$ on the $z=0$ plane and radius $t$ ($(P,t)=(p,q,t)\in \mathbf R^{2}\times(0,\infty)$ ).

\begin{theorem}\label{thm2}  For $(x,y,z)\in\mathbf R^z$, where $z\ne0$ we have
	\begin{eqnarray}
		f(x,y,z)=
		\lim_{n\to \infty}\left((n+1)\times\right.\,\,\,\,\,\,\,\label{1.6}\\\left.
		\left((2n+1) Mf(x,y,\mid{z}\mid)+(\sgn z)(\frac{2}{3}n+1)a_{01}(x,y,\mid{z}\mid)\right)+\sum_{i=0}^{n}\int_{0}^{\mid{z}\mid}{\mid{z}\mid}^{2i-1}\nonumber\right.\\\left.\times
		\left(C_{2n,2i}
		(\frac{u}{\mid{z}\mid})\Delta^{i}(Mf(x,y,u)) +(\sgn z)C_{2n+1,2i}
		(\frac{u}{\mid{z}\mid})\Delta^{i}(a_{01}(x,y,u))\right)du\right)
		\nonumber		\end{eqnarray}
	here sgn is the signum function,  $\Delta$  is the Laplace operator in two dimensions  with respect to $x$ and $y$ ($\Delta^{0}(\cdot)=(\cdot)$).	
\end{theorem}

\noindent In Section 3, we consider the spherical harmonic expansion of the restriction of $f$ onto $S(x, y, z)$. The limit is convergent because the expression in the brackets in (\ref{1.6}) is the $n$-th partial sum of the  expansion (see \ref{3.2}  below).

\noindent  From Theorem 1.2 it follows that for reconstructing $f(x,y,z)$ at any point $(x,y,z)$ we use the local data  of $Mf$ and $a_1$. 

\begin{corollary} Let $f\in{ C}^{\infty}(\mathbf R^3)$ to restore  $f(x,y,z)$ for any $(x,y,z)\in\mathbf R^3$ we use the values of $Mf$ and $a_{01}$ over spheres with $(p,q,t)\in U\times (0,\mid{z}\mid]$ where $U\subset\{z=0\}$ is a neighborhood of  $(x,y,0)\in\{z=0\}$.
\end{corollary}

\noindent Also, note that using (\ref{1.6})	one can find an inversion formula of SRT for functions supported on one side of a plane, not necessarily with compact support (see \cite{Ara231}). Indeed, let $f_c$ be a smooth function defined on $\mathbf R^3$, supported in  halfspace $\{z>0\}$. We consider the following even function (with respect to the $\{z=0\}$ plane)

\begin{equation}\label{1.7}	\begin{cases}
		f(\textbf{x})=f_c(\textbf{x}),\,\,\,\text{for}\,\,\,\,  \textbf{x}\in\{z\geq0\}\\		
		f(\textbf{x})=f_c(-\textbf{x}),\,\,\,\,\text{for}\,\,\,  \textbf{x}\in\{z\leq0\}.
	\end{cases}	
\end{equation}
\noindent It is clear that for the Fourier first cosine coefficient of the restriction of  $f$ onto the sphere $S(x,y,t)$, we have $a_{01}(x,y,t)\equiv0$. Hence $f$ (as well as $f_c$) can be restored from its SRT over spheres centered on the $\{z=0\}$ plane.

\noindent An illustration of the algorithm obtained by (\ref{1.6}) is provided in Section 8.

\section{Injectivity of the two-data spherical Radon transform}\label{sec2}

\noindent \begin{proof}{(Theorem 1.1)}	 We will show: if
	\begin{equation}\label{2.3}
		a_{00}(p,q,t)=0 \,\,\,\,\text{and}\,\,\,\,
		a_{01}(p,q,t)=0
	\end{equation}
	for a given function $f\in {\textit{C}}^{1}(\mathbf R^2)$ (for all $(p,q,0)\in D$ and $t>0$) then $f\equiv0$.
	
	\noindent Derivating (\ref{2.3})  with respect to $p$ and $q$ we get (using (\ref{1.3.1}))
	\begin{equation}\label{2.5}
		\int_{S(p,q,t)} f'_x\,d\omega=0\,\,\, \, \text{and}\,\,\,\,\,
		\int_{S(p,q,t)} f'_y\,d\omega=0.\end{equation}
	
	\noindent	Using (\ref{2.5}) and	Ostrogradsky's theorem we get ($P_{1,1}(\cos\theta)=\sin\theta$)
	\begin{eqnarray}
		t^2\int_{S(p,q,t)} f \sin\theta\cos\varphi d\omega=\int_{B(p,q,t)} f'_x dxdydz=
		\int_0^t u^2du\int_{S(p,q,u)} f'_x d\omega=0\nonumber
	\end{eqnarray}
	here $B(p,q,t)$ is the ball of radius $t$ centered at $P=(p,q,0)$. Hence 
	
	\noindent $a_{11}(p,q,t)=0$ for all $(p,q,0)\in D$ and $t>0$.
	In the same way, we get $b_{11}(p,q,t)=0$ for all $(p,q,0)\in D$ and $t>0$.

	\noindent Now, we assume for all $(p,q,0)\in D$ and $t>0$  ($n\geq 1$)
	\begin{equation}\label{2.14}\begin{cases}
			a_{mk}(p,q,t)=0\,\,\, \text{for}\,\,\, 0\leq m\leq k, \,\,\,\,k\leq n\\
			b_{mk}(p,q,t)=0\,\,\, \text{for}\,\,\, 0\leq m\leq k, \,\,\,\,k\leq n
		\end{cases}
	\end{equation}
	and we are going to show (\ref{2.14}) for $k=n+1$ (the mathematical induction method).  
	
	\noindent Using recurrence properties for the associated Legendre polynomials (see \cite{Pou}), for $k=n+1$ and $0< m\leq n+1$ we have 	
	\begin{eqnarray}
		\frac{2\pi}{2n+3}	a_{m(n+1)}(p,q,t)=(-1)^m\int_{S(p,q,t)} f\,P_{n+1}^m(\cos\theta)\,\cos(m\varphi)\,d\omega=\nonumber\\
		\frac{(-1)^{m+1}}{2m}\int_{S(p,q,t)} f\,\sin\theta[P_{n}^{m+1}+(n+m+1)(n+m)P_{n}^{m-1}]\,\cos(m\varphi)d\omega.\label{2.17}\end{eqnarray}	
	
	\noindent  Using Ostrogradsky's theorem for the second term of (\ref{2.17}) we get	
	\begin{eqnarray}
		II=\int_{S(p,q,t)} f\sin\theta P_{n}^{m-1}(\cos\theta)[\cos(m-1)\varphi\cos\varphi-\sin(m-1)\varphi\sin\varphi]d\omega=\label{2.18}\\							\int_{B(p,q,t)}f[(P_{n}^{m-1}(\cos\theta)\cos(m-1)\varphi)'_x -(P_{n}^{m-1}(\cos\theta)\sin(m-1)\varphi)'_y]\,dxdydz.\,\,\,\,\label{2.18}\nonumber\end{eqnarray}	
	
	\noindent 	Representating $\sin (m-1)\varphi$ and  $\cos (m-1)\varphi$ in terms of the powers of $\sin\varphi$ and $\cos \varphi$ respectively one can calculate that for any  integer $0\leq m-1$ (see \cite{Ara232} eq. (2.16))
	\begin{equation}
		(P_{n}^{m-1}\,\cos(m-1)\varphi)'_x-(P_{n}^{m-1}\,\sin(m-1)\varphi)'_y=0.
		\label{2.21}
	\end{equation}
	
	\noindent Hence $II=0$ for any $0\leq m-1\leq n$. For the first term of (\ref{2.17}), using the recurrence properties for the associated Legendre polynomials (see \cite{Pou}), we have 	
	
	\begin{eqnarray}
		I=\int_{S(p,q,t)} f\sin\theta P_{n}^{m+1}(\cos\theta)\cos(m\varphi)d\omega=\int_{S(p,q,t)} f\sin\theta\cos(m\varphi)\times\nonumber\\ (P_{n-2}^{m+1}+(n+m-1)(n+m-2)P_{n-2}^{m-1}-(n-m)(n-m+1)P_{n}^{m-1})
		d\omega\nonumber\\
		=\int_{S(p,q,t)} f\sin\theta P_{n-2}^{m+1}\cos(m\varphi)d\omega.\label{2.23}\end{eqnarray}	
	
	\noindent  Using the same  recurrence properties one can continue the process until $n-2< m+1$ and taking into account that $P_{n}^{m}=0$ for $n<m$ we obtain $I=0$.
	
	\noindent Hence, for all $(p,q)\in D$, $t>0$ and $0<m\leq n+1$  we obtain 
	\begin{equation}\label{2.25}
		a_{m(n+1)}(p,q,t)=0.
	\end{equation} 
	In a similar way, for $0<m\leq n+1$, one can prove that 	$b_{m(n+1)}(p,q,t)=0$.

	\noindent	Now, we prove (\ref{2.14}) for $k=n+1$ and $m=0$. It is known that $P_{n}^{0}=
	P_{n}(\cos\theta)=\sum_{k=o}^{[n/2]} c_k(\cos\theta)^{n-2k}$,
	where $c_k$ is a constant (see \cite{Pou}). Hence, it is sufficient  
	to prove: from (for all $(p,q,0)\in D$ and $t>0$)
	\begin{equation}\label{2.28}         
		\int_{S(p,q,t)} f \cdot(\cos\theta)^{k} d\omega=0 \,\, \texttt{for}\,\,0\leq k\leq n,\,\,\,\, n\geq 1
	\end{equation}
	follows (\ref{2.28}) for $k=n+1$.
	Using Ostrogradsky's theorem we get	
	\begin{equation} \label{2.29}      
		\int_{S(p,q,t)} f \cdot(\cos\theta)^{n+1} d\omega=\int_{S(p,q,t)} f \cdot(\frac{z}{t})^n\cos\theta d\omega=	\int_{B(p,q,t)} f'_z\cdot(\frac{z}{t})^n \,dxdydz.\end{equation}	
	
	\noindent It follows from (\ref{2.28}) that       
	\begin{equation}\label{2.30}  \int_{S(p,q,t)} f'_x \cdot(\cos\theta)^{n-1} d\omega=0,\,\,\,\,\,\,\int_{S(p,q,t)} f'_y \cdot(\cos\theta)^{n-1}d\omega=0.
	\end{equation} 
	First applying Ostrogradsky's theorem to (\ref{2.30}), then differentiating it with respect to $p$ and $q$ respectively we get
	\begin{equation}\label{2.32}         
		\int_{S(p,q,t)} f'_x \cdot(\cos\theta)^{n-1} \sin\theta\cos\varphi d\omega=0,\,\,\, \,\,\, \int_{S(p,q,t)} f'_y \cdot(\cos\theta)^{n-1}\sin\theta\sin\varphi d\omega=0.
	\end{equation}		
	Now, differentiating $	\int_{S(p,q,t)} f \cdot(\cos\theta)^{n-1} d\omega=0 $ with respect to $t$ and using (\ref{2.32}) and (\ref{2.29}) we obtain  (\ref{2.28}) for $k=n+1$.	
	Thus for any $(p,q,0)\in D$ and $t>0$ we obtain
	\begin{equation}  \label{2.27}       
		a_{m(n+1)}(p,q,t)=b_{m(n+1)}(p,q,t)=0
	\end{equation}
	for all $n\geq0$ and $0\leq m\leq n$,  hence $f\equiv0$. \end{proof}

\noindent We have proved Theorem 1 directly.

\noindent	The problem of reconstructing a function $f$ from the spherical means is known
to be equivalent (via the Kirchoff-Poisson formula) to the following boundary value problems for the wave equation (\cite{AgQu}-\cite{AmKu}, \cite{KuKu}, \cite{St}) 
\begin{equation}\label{2.35}\begin{cases}
		u(\textbf{x},t)''_{tt}-\Delta_\textbf{x} u(\textbf{x},t)=0,\,\,\,\,  (\textbf{x},t)\in\mathbf R^{d}\times(0,\infty) \\
		u(\textbf{x},t)=f(\textbf{x}),\,\, \,\,\,\,	u(\textbf{x},0)'_t=0,
	\end{cases}
\end{equation}
where $u(\textbf{x},t)$ is the induced pressure wave,
$f$ is the initial pressure and the problem is to find the initial value $f(\textbf{x}) = u(\textbf{x}; 0)$ from the Dirichlet
data (the solution
of the wave equation (\ref{2.35}) on $L\times(0,\infty)$). 

\noindent	The assumption that the support of $f$ is compact becomes decisive in solving both the integral equation (\ref{1.0}) and the wave equation (\ref{2.35}) (\cite{KuKu}). 

\noindent	In the space of compactly supported functions one can prove Theorem 1.1 using  known general results in PDE. Now, we are trying to write the boundary value conditions for the wave equation (\ref{2.35}) equivalent to the problem we are studying in the article.

\noindent One can prove that the functions 	$a_{00}(p,q,t)$ and $a_{01}(p,q,t)$ are related by the following relation which follows from the Stokes divergence formula. 

\begin{lemma}\label{0} 
	\begin{equation}\label{2.36}
		t^2\dfrac{\partial}{\partial n}(3 a_{00}(p,q,t))	=
		\dfrac{\partial}{\partial t}(t^2 a_{01}(p,q,t)). \,\,\,\,\, (P,t)\in D\times (0,\infty),
	\end{equation}
	where $\dfrac{\partial}{\partial n}(\cdot)$ is the normal derivative at $P\in D$ with respect
	to the spatial variable.
\end{lemma}

\noindent Hence, it follows from Lemma 1 that the data $a_{01}(p,q,t)$ expresses
the normal derivative ${\partial_n}(a_{00}(p,q,t))$ and therefore is known when $P\in D$. Thus, this means that both of the Cauchy (Dirichlet-Neumann) data of the solution of the wave equation $u: \mathbf R^3 \times[0,\infty) \to \mathbf R^1$ along $D$ with respect
to the spatial variable \begin{equation}\label{2.37}
	u(\textbf{x}; t)\,\,\, \text{ and}\,\,\, \partial_n(u(\textbf{x},t)), \,\,\, (\textbf{x},t)\in D\times(0,\infty) \end{equation}
are given. Now, in terms of PDE the problem considered in the article is to study the problem of  reconstructing $f$ in wave equation (\ref{2.35})
from Cauchy (Dirichlet-Neumann) conditions (\ref{2.37}) on $D$. 
The uniqueness is followed from the Holmgren-John theorem (see \cite{Kur}). 

\section{Reconstruction formula using the Fourier coefficients}\label{3}

\noindent  Let 

\noindent	$(x,y,z)\in \mathbf R^3$, be a point with $z>0$. On $S(x,y,z)$ (the sphere with center $P=(x,y,0)$ and radius $t=z$) for the polar coordinates of $(x,y,z)$ with respect to $(x,y,0)$ we have 
$$(x,y,z)=(t,\theta,\varphi)=(z, 0,\varphi).$$
Hence, using (\ref{1.4}), taking into account that $P_n(\cos0)=P_n(1)=1$ and $P_{nm}(\cos0)=P_{nm}(1)=0$ for $m\geq1$, we get
\begin{equation}\label{3.1}
	f(x,y,z)=f_{(x,y,z)}(0,\varphi)=Mf(x,y,z)+\sum_{n=1}^{\infty}a_{0n}(x,y,z).\end{equation}

\noindent For $(x,y,z)\in\mathbf R^3$, with $z<0$, on $S(x,y,\mid{z}\mid)$ 
for the polar coordinates of $(x,y,z)$ with respect to $(x,y,0)$ we have 
$$(x,y,z)=(z, \pi,\varphi).$$
Hence, using (\ref{1.4}) and taking into account that 	 $\,P_n(-1)=(-1)^nP_n(1)=(-1)^n$ and $P_{nm}(-1)=0$ for $m\geq1$, we get
\begin{eqnarray}\label{3.2}
	f(x,y,z)=f_{(x,y,z)}(\pi,\varphi)=Mf(x,y,\mid{z}\mid)+\\\sum_{n=1}^{\infty}a_{0n}(x,y,\mid{z}\mid)P_n(-1)=
	Mf(x,y,\mid{z}\mid)+\sum_{n=1}^{\infty}(-1)^na_{0n}(x,y,\mid{z}\mid)\nonumber.\end{eqnarray}
Thus for $(x,y,z)\in \mathbf R^3$, where $z\ne0$, we get
\begin{equation}\label{3.3}	f(x,y,z)=
	Mf(x,y,\mid{z}\mid)+\sum_{n=1}^{\infty}(\sgn{z})^na_{0n}(x,y,\mid{z}\mid)\nonumber.\end{equation}

\noindent  Thus to reconstruct $f$  we need to find $a_{0n}$ for $n>2$. Note that in  (\cite{Ara231}) to find an inversion formula of SRT for functions supported on one side of a plane we needed to find $a_{0(2n)}$ for $n>1$.

\section{{The derivation of the Fourier coefficients}}\label{sec4}

\noindent For $(x,y,z)\in \mathbf R^3$, we consider the bundle of spheres $S(P,t)$ where $P=(p,q,0)\in \{z=0\}$ containing $(x,y,z)$. 	
Any sphere from the bundle is determined by its center $(p,q,0)\in \{z=0\}$. Let $(t(p,q),\theta(p,q),\varphi(p,q))$ be the polar coordinates of $(x,y,z)$ with respect to $(p,q,0)\in \{z=0\}$, in short we will write $(t(p,q),\theta(p,q),\varphi(p,q))=(x,y,z)$.

\noindent The consistency method is based on the following statement. The restrictions  

\noindent	$f_{(p,q,t)}(\theta,\varphi)$ are consistent. We mean the following: for a fixed 
$(x,y,z)\in \mathbf R^3$, we have

\begin{equation}\label{4.0}
	f(x,y,z)=f_{(p,q,t)}(\theta,\varphi) \,\,\,\,\emph{for all}\,\,\,\,\, (p,q,t,\theta,\varphi))=(x,y,z)\end{equation}
hence (\ref{4.0}) does not depend on $p$ and $q$ (the consistency condition).

\noindent In this article we derive the unknown Fourier coefficients of the restriction of  $f$ onto the sphere $S(p,q,t)$ from the  consistency condition of $f_{(p,q,t)}(\theta,\varphi)$ (see \cite{Ara10},\cite{Ara19}, \cite{Ara231},\cite{Ara232}).	

\noindent For a fixed $(x,y,z)$ differentiating (\ref{1.3.1}) with respect to $p$ we get
\begin{equation}\label{4.2}\begin{cases}
		0=1+t'_p \sin\theta\cos\varphi+t\cos\theta\cos\varphi\,\theta'_p-t\sin\theta\sin\varphi\,\varphi'_p\\
		0=t'_p \sin\theta\sin\varphi+t\cos\theta\sin\varphi\,\theta'_p+t\sin\theta\cos\varphi\,\varphi'_p\\
		0=t'_p \cos\theta-t\sin\theta\,\theta'_p.\end{cases}
\end{equation}
From (\ref{4.2}), we obtain
\begin{equation}\label{4.3}
	t'_p=-\cos\varphi\sin\theta, \quad \varphi'_p=\frac{\sin\varphi}{t\sin\theta},\quad \theta'_p=-\frac{\cos\varphi\cos\theta}{t}.
\end{equation}

\noindent Similarly, differentiating (\ref{1.3.1}) with respect to $q$, we get
\begin{equation}\label{4.4}
	t'_q=-\sin\varphi\sin\theta, \quad \varphi'_q=-\frac{\cos\varphi}{t\sin\theta},\quad \theta'_q=-\frac{\sin\varphi\cos\theta}{t}.
\end{equation}

\noindent Differentiating  (\ref{4.0}) we obtain the following Lemma.
\begin{lemma}\label{1} Using the  consistency condition we obtain the following equations for the restrictions $f_{(p,q,t)}(\theta,\varphi)$ 
	\begin{equation}\label{4.5}
		(f(x,y,z))'_p=
		f'_p-f'_t\cos\varphi\sin\theta+f'_\varphi\frac{\sin\varphi}{t\sin\theta}-
		f'_\theta\frac{\cos\varphi\cos\theta}{t}=0
	\end{equation}
	and
	\begin{equation}\label{4.6}
		(f(x,y,z))'_q=
		f'_q-f'_t\sin\varphi\sin\theta-f'_\varphi\frac{\cos\varphi}{t\sin\theta}-f'_\theta\frac{\sin\varphi\cos\theta}{t}=0.
	\end{equation}
\end{lemma}

\noindent In order to get differential equations for  $a_{0n}(p,q,t)$ and 
$a_{1n}(p,q,t)$ we
multiply (\ref{4.5}) by $P_n(\cos\theta)$, for $n\geq0$, and integrate over ${\mathbf S^{2}}$. We obtain

\begin{eqnarray}
	\int_{\mathbf S^{2}}f'_p\,P_n(\cos\theta)d\omega-\int_{\mathbf S^{2}}f'_t\cos\varphi\sin\theta\,P_n(\cos\theta)d\omega+\label{4.7}\\
	\frac{1}{t}\int_{\mathbf S^{2}}f'_\varphi\,{\sin\varphi}\,\frac{P_n(\cos\theta)}{\sin\theta}d\omega -\frac{1}{t}\int_{\mathbf S^{2}}f'_\theta\,\cos\varphi\,\cos{\theta}\,P_n(\cos\theta)d\omega=0.\nonumber\end{eqnarray}

Considering (\ref{4.7}) term by term, using Integration by parts and known recurrence relations for the associated Legendre polynomials (see (\cite{Pou}) ), for $n\geq0$ we get (detailed calculation can be found in (\cite{Ara231}))
\begin{eqnarray}
	\frac{(n+1)(n+2)}{2n+3}(a_{1(n+1)})'_t+\frac{(n+1)(n+2)^2}{t(2n+3)}a_{1(n+1)}-\frac{n(n-1)}{2n-1}(a_{1(n-1)})'_t+\label{4.14}\\
	\frac{n(n-1)^2}{t(2n-1)}a_{1(n-1)}-2(a_{0n})'_p=0\nonumber\end{eqnarray}
here we assume $a_{1(-1)}\equiv a_{1(0)}\equiv0$ and $a_{00}(p,q,t)=Mf(p,q,t)$, since $P_0(\cos\theta)=1$.

\noindent In the same way, multiplying (\ref{4.6}) by $P_n(\cos\theta)$ and integrating over ${\mathbf S^{2}}$ for $n\geq0$ we derive 

\begin{eqnarray}
	\frac{(n+1)(n+2)}{2n+3}(b_{1(n+1)})'_t+\frac{(n+1)(n+2)^2}{t(2n+3)}b_{1(n+1)}-\frac{n(n-1)}{2n-1}(b_{1(n-1)})'_t+\label{4.16}\\
	\frac{n(n-1)^2}{t(2n-1)}
	b_{1(n-1)}-2(a_{0n})'_q=0\nonumber\end{eqnarray}
here we assume $b_{1(-1)}\equiv b_{1(0)}\equiv0$.

\noindent We need to derive the next equation for the
unknown Fourier coefficients. Multiplying (\ref{4.5}) and (\ref{4.6}) by $\cos\varphi\,P_{n1}(\cos\theta)$ and $\sin\varphi\,P_{n1}(\cos\theta)$ ($n\geq1$), respectively, and integrating their addition over ${\mathbf S^{2}}$,  we obtain 

\begin{eqnarray}
	\int_{\mathbf S^{2}}(f'_p\,\cos\varphi\,P_{n1}(\cos\theta)+f'_q\,\sin\varphi\,P_{n1}(\cos\theta))d\omega-\int_{\mathbf S^{2}}f'_t\sin\theta\,P_{n1}(\cos\theta)d\omega\nonumber\\
	-\frac{1}{t}\int_{\mathbf S^{2}}f'_\theta\,\cos{\theta}\,P_{n1}(\cos\theta)d\omega=0.\label{4.17}\end{eqnarray}

Considering (\ref{4.17}) term by term, using Integration by parts and known recurrence relations for the associated Legendre polynomials (see (\cite{Pou})), for $n\geq1$ we get 

\begin{eqnarray}	\label{4.21}
	\frac{2}{2n+3}(a_{0(n+1)})'_t+\frac{2(n+2)}{t(2n+3)}a_{0(n+1)}-\frac{2}{2n-1}(a_{0(n-1)})'_t+\frac{2(n-1)}{t(2n-1)}a_{0(n-1)}+\nonumber\\
	(a_{1n})'_p+(b_{1n})'_q=0.
\end{eqnarray}

\noindent It follows from (\ref{4.21}) that to calculate $a_{0n}$ we need to know $(a_{1n})'_p+(b_{1n})'_q$. Differentiating by $p$ and $q$ respectively and taking the sum of the corresponding equations of (\ref{4.14}) and (\ref{4.16}) for $n\geq0$ we obtain

\begin{eqnarray}
	\frac{(n+1)(n+2)}{2n+3}((a_{1(n+1)})'_p+(b_{1(n+1)})'_q)'_t+\frac{(n+1)(n+2)^2}{t(2n+3)}((a_{1(n+1)})'_p+(b_{1(n+1)})'_q)-\nonumber\\\frac{n(n-1)}{2n-1}((a_{1(n-1)})'_p+(b_{1(n-1)})'_q)'_t+\label{4.22}\\
	\frac{n(n-1)^2}{t(2n-1)}
	((a_{1(n-1)})'_p+(b_{1(n-1)})'_q)-2\Delta(a_{0n})=0\nonumber
\end{eqnarray}
here we assume $b_{1(-1)}\equiv b_{1(0)}\equiv a_{1(-1)}\equiv a_{1(0)}\equiv0$. 	

\noindent
Finally, to find $a_{0(n)}$ for any $n>1$ we apply the following algorithm.  

\noindent The algorithm: First, we find $(a_{11})'_p+(b_{11})'_q$ from (\ref{4.22}) (written for $n=0$) and by substituting it into (\ref{4.21}) (written for $n=1$), we can calculate $a_{02}$. Then, we find $(a_{12})'_p+(b_{12})'_q$ from (\ref{4.22}) (written for $n=1$) and by substituting it into (\ref{4.21}) (written for $n=2$), we derive $a_{03}$ (note, that $a_{01}$ is known). Next, substituting $a_{02}$, $(a_{13})'_p+(b_{13})'_q$ (obtained from (\ref{4.22}) (written for $n=2$)) into (\ref{4.21} (written for $n=3$), we calculate $a_{04}$ and so on. Now, we are going to apply this algorithm.

\section{{A representation for the Fourier coefficients}}\label{sec5}

\noindent From (\ref{1.4.1}) we have the following boundary conditions (for  $n\geq1, 1\leq m\leq n$)
\begin{equation}\label{5.1}
	a_{0n}(p,q,0)=0,  \,\, a_{mn}(p,q,0)=0, \,\,b_{mn}(p,q,0)=0.\end{equation}

\noindent For $a_{0(2k+1)}$, ($k\geq1$) we get (the  solution of (\ref{4.21}) with the boundary conditions (\ref{5.1}))
\begin{eqnarray}
	a_{0(2k+1)}(p,q,t)=
	\frac{(4k+3)}{4k-1} a_{0(2k-1)}+\dfrac{1}{t^{2k+2}}\int_{0}^{t}(4k+3)u^{2k+1}\times \nonumber\\\left(-\dfrac{4k+1}{4k-1}a_{0(2k-1)} -\frac{u}{2}((a_{1 (2k)}(p,q,u))'_p+(b_{1 (2k)}(p,q,u))'_q)\right)du.\label{5.2}
\end{eqnarray}

\noindent Also, for $(a_{1(2k)}(p,q,u))'_p+(b_{1(2k)}(p,q,u))'_q$, ($k\geq1$) we get (the  solution of (\ref{4.22}) with the boundary conditions (\ref{5.1}))
\begin{eqnarray}
	(a_{1(2k)}(p,q,t))'_p+(b_{1(2k)}(p,q,t))'_q=\frac{(4k+1)(2k-1)(2k-2)}{2k(4k-3)(2k+1)}\times\label{5.2.1}\\((a_{1(2k-2)}(p,q,t))'_p+(b_{1(2k-2)}(p,q,t))'_q)+
	\frac{1}{t^{2k+1}}\int_{0}^{t}\frac{(4k+1)u^{(2k)}}{2k(2k+1)}\times\nonumber\\
	\left(-\dfrac{(2k-1)(2k-2)(4k-1)}{(4k-3)} ((a_{1(2k-2)}(p,q,u))'_p+(b_{1(2k-2)}(p,q,u))'_q)+\right.\nonumber\\\left.2u\Delta(a_{0(2k-1
		)}(p,q,u))\right)d\,u.\nonumber\end{eqnarray}

\begin{theorem}\label{1} The following representations are valid. For $k>0$ we have
	
	\begin{eqnarray}\label{5.3}
		a_{0(2k-1)}(p,q,t)=\dfrac{(4k-1)}{3} a_{01}(p,q,t)+\\ \int_{0}^{t}\sum_{i=0}^{k-1}t^{2i-1}\,\sum_{m=1}^{k+i-1}c_m(2k-1,2i)\,(u/t)^{2m+1}\Delta^{i}(a_{01}(p,q,u))\, du,\nonumber
	\end{eqnarray}
	and for $k>1$	                 
	
	\begin{eqnarray}\label{5.4}
		(a_{1(2k-2)}(p,q,u))'_p+(b_{1(2k-2)}(p,q,u))'_q=\\
		\int_{0}^{t}\sum_{i=1}^{k-1}t^{2(i-1)}\,\sum_{m=1}^{k+i-2}s_m(2(k-1),2(i-1))\,(u/t)^{2m+1}(\Delta^{i}(a_{01}(p,q,u)))\, du.\nonumber
	\end{eqnarray}
\end{theorem}
\textbf{	Proof.}  From (\ref{4.22}) we have

\begin{equation}\label{5.5}
	(a_{12}(p,q,t))'_p+(b_{12}(p,q,t))'_q=\int_{0}^{t} \dfrac{5}{3}(\frac{u}{t})^3\Delta(a_{01}(p,q,u))) du,
\end{equation}	
hence (\ref{5.3}) for $k=1$ and for $k=2$ (\ref{5.4}) are true.
Suppose (\ref{5.4}) is true for some $n = k$. Show that  it is true for $n = k + 1$ as well (the mathematical induction method).                                                                                                                  	

\noindent After substituting $(a_{1(2k-2)})'_p+(b_{1(2k-2)})'_q$ and $a_{0(2k-1)}$ from (\ref{5.4}) and (\ref{5.3}), respectively, into (\ref{5.2}) and simplifying the expression, we change the order of summation, then the order of integration, and after grouping we get
\begin{eqnarray}
	\frac{2k(2k+1)}{4k+1}	[(a_{1(2k)}(p,q,t))'_p+(b_{1(2k)}(p,q,t))'_q]=\nonumber\\
	\int_{0}^{t}\frac{8k-2}{3}\Delta(a_{01}(p,q,u))(\frac{u}{t})^{2k+1}du\nonumber\\+\sum_{i=1}^{k}t^{2(i-1)}\int_{0}^{t}\left(\left(\sum_{m=1}^{k+i-2}
	\frac{(2k-1)(2k-2)s_m(2(k-1),2(i-1))}{4k-3}(\frac{u}{t})^{2m+1}+\right.\right.\nonumber\\\left.\left.                                                                                                                                                                         \sum_{m=1}^{k+i-2}\left(\frac{c_m(2k-1,2(i-1))}{k+i-m-1}-\frac{(2k-1)(k-1)(4k-1)s_m(2(k-1),2(i-1))}{(4k-3)(k+i-m-1)}\right)\nonumber\right.\right.\\\left.\left.\times[(\frac{u}{t})^{2m+1}-(\frac{u}{t})^{2k+2i-1}]\right)(\Delta^{i}(a_{01}(p,q,u)))\right)du	\,\,\,\,\,\,\,\,\,\,\,\,\,	\label{5.6}	\end{eqnarray}
hence (\ref{5.2}) is true for every $k$.	

\noindent After substituting $(a_{1(2k-2)})'_p+(b_{1(2k-2)})'_q$ and $a_{0(2k-1)}$ from (\ref{5.4}) and (\ref{5.3}), respectively, into (\ref{5.2}), and simplifying the expression, we change the order of summation, then the order of integration, and after grouping we get

\begin{eqnarray}
	\frac{a_{0(2k+1)}(p,q,t)}{4k+3}=
	\frac{a_{01}(p,q,t)}{3}+\int_{0}^{t}t^{-1}a_{01}(p,q,u)\left(\frac{-(4k+1)}{3}(\frac{u}{t})^{2k+1}+\right.\nonumber\\
	\left.\sum_{m=1}^{k-1}\left(\frac{c_m(2k-1,0)}{4k-1}(\frac{u}{t})^{2m+1} -\frac{(4k+1)c_m(2k-1,0)}{(8k-2)(k-m)}
	[(\frac{u}{t})^{2m+1}-(\frac{u}{t})^{2k+1}]\right)\right)du
	\nonumber\\
	+\sum_{i=1}^{k}t^{2i-1}\int_{0}^{t}(\Delta^{i}(a_{01}(p,q,u)))\sum_{m=1}^{k+i-1}\left(\frac{c_m(2k-1,2i)}{(4k-1)}(\frac{u}{t})^{2m+1}-\right.\nonumber\\\left.\left(\frac{(4k+1)c_m(2k-1,2i)}{(8k-2)(k+i-m)}+\frac{s_m(2k,2(i-1))}{4(k+i-m)}\right)[(\frac{u}{t})^{2m+1}-(\frac{u}{t})^{2k+2i+1}]\right)du.\,\,\,\,\,\,\,\label{5.7}	
\end{eqnarray}	
Theorem 3 is proved.

\noindent For $a_{0(2k+2)}$, ($k\geq0$) we get (the  solution of (\ref{4.21}) with the boundary conditions (\ref{5.1}))

\begin{eqnarray}
	a_{0(2k+2)}(p,q,t)=
	\frac{(4k+5)}{4k+1} a_{0(2k)}(p,q,t)+\dfrac{1}{t^{2k+3}}\int_{0}^{t}(4k+5)u^{2k+2}\times \label{5.8}\\\left(-\dfrac{4k+3}{4k+1}a_{0(2k)} -\frac{u}{2}((a_{1 (2k+1)}(p,q,u))'_p+(b_{1 (2k+1)}(p,q,u))'_q)\right)du.\nonumber
\end{eqnarray}

\noindent Also, for $(a_{1(2k+1)}(p,q,u))'_p+(b_{1(2k+1)}(p,q,u))'_q$, ($k\geq1$) we get (the  solution of (\ref{4.22}) with the boundary conditions (\ref{5.1}))

\begin{eqnarray}
	\frac{(2k+1)(2k+2)}{4k+3}	(a_{1(2k
		+1)}(p,q,t))'_p+(b_{1(2k+1)}(p,q,t))'_q=\label{5.8.1}\\\frac{2k(2k-1)}{(4k-1)}((a_{1(2k-1)}(p,q,t))'_p+(b_{1(2k-1)}(p,q,t))'_q)+
	\frac{1}{t^{2k+2}}\int_{0}^{t}u^{(2k+1)}\nonumber\\\times
	\left(-\dfrac{(2k-1)(2k)(4k+1)}{(4k-1)} ((a_{1(2k-1)}(p,q,u))'_p+(b_{1(2k-1)}(p,q,u))'_q)+\right.\nonumber\\\left.2u\Delta(a_{0(2k)}(p,q,u))\right)d\,u;\nonumber\end{eqnarray}

\begin{theorem}\label{2} The following representations are valid. For $k>0$ we have		
	\begin{eqnarray}\label{5.9}
		a_{0(2k)}(p,q,t)=(4k+1) Mf(p,q,t)+\\ \int_{0}^{t}\sum_{i=0}^{k}t^{2i-1}\,\sum_{m=1}^{k+i}c_m(2k,2i)\,(u/t)^{2m}\Delta^{i}(Mf(p,q,u))\, du,\nonumber
	\end{eqnarray}
	and 
	\begin{eqnarray}\label{5.10}
		(a_{1(2k-1)}(p,q,u))^\prime_p+(b_{1(2k-1)}(p,q,u))^\prime_q=\\\int_{0}^{t}\sum_{i=1}^{k}t^{2(i-1)}\,\sum_{m=1}^{k+i-1}s_m(2k-1,2(i-1))\,(u/t)^{2m}(\Delta^{i}(Mf(p,q,u)))\, du,\nonumber
	\end{eqnarray}
\end{theorem}
\textbf{Proof.}  For $k=1$ (\ref{5.10}) and (\ref{5.9}) are true. From (\ref{4.22}) we have

\begin{equation}\label{5.11}
	(a_{11}(p,q,u))^\prime_p+(b_{11}(p,q,u))^\prime_q=\int_{0}^{t} 3(\frac{u}{t})^2\Delta(Mf(p,q,u)) du
\end{equation}
and from (\ref{4.21}) we have
\begin{eqnarray}
	a_{02}(p,q,t)= 5Mf(p,q,t)+\qquad\label{5.12}\\ \int_{0}^{t}\left(t^{-1}(-15(\frac{u}{t})^2)Mf(p,q,u)+t(\frac{15}{4})[(\frac{u}{t})^4-(\frac{u}{t})^2]\Delta(Mf(p,q,u))\right) du.\nonumber\end{eqnarray}
Thus (\ref{5.10}) and (\ref{5.9}) are true for $k=1$

Suppose (\ref{5.9})-(\ref{5.10}) are true for some $n = k$. 

\noindent After substituting $(a_{1(2k-1)})'_p+(b_{1(2k-1)})'_q$ and $a_{0(2k)}$ from (\ref{5.9}) and (\ref{5.10}), into (\ref{5.8.1}), and simplifying the expression, we change the order of summation, then the order of integration, and after grouping them, we get

\begin{eqnarray}
	\frac{(2k+1)(2k+2)}{4k+3}	(a_{1(2k+1)})'_p+(b_{1(2k+1)})'_q=\int_{0}^{t}\left((8k+2)(\frac{u}{t})^{2k+2}\times\right.\nonumber\\\left.
	(\Delta(Mf(p,q,u)))+\sum_{i=1}^{k+1}t^{2(i-1)}\left(\sum_{m=1}^{k+i-1}\dfrac{2k(2k-1)s_m(2k-1,2(i-1))}{4k-1}(\frac{u}{t})^{2m}\right.\right.\nonumber\\\left.\left.+\sum_{m=1}^{k+i-1}\left(\dfrac{c_m(2k,2(i-1))}{k+i-m}-\dfrac{k(2k-1)(4k+1)s_m(2k-1,2(i-1))}{(4k-1)(k+i-m)}\right)\times\nonumber\right.\right.\\\left.\left.[(\frac{u}{t})^{2m}-(\frac{u}{t})^{2k+2i}]\right)(\Delta^{i}(Mf(p,q,u)))\right)du	\,\,\,\,\,\,	\label{5.13}	\end{eqnarray}

\noindent After substituting $(a_{1(2k-1)})'_p+(b_{1(2k-1)})'_q$ and $a_{0(2k)}$ from (\ref{5.9}) and (\ref{5.10}), into (\ref{5.8}),and simplifying the expression, we change the order of summation, then the order of integration, and after grouping them, we get

\begin{eqnarray}
	\dfrac{a_{0(2(k+1))}(p,q,t)}{4k+5}=Mf(p,q,t)-\int_{0}^{t}t^{-1}Mf(p,q,u)\left((4k+3)(\frac{u}{t})^{2k+2}\right.\nonumber\\
	\left.-\sum_{m=1}^{k}\dfrac{c_m(2k,0)}{4k+1}(\frac{u}{t})^{2m} +\sum_{m=1}^{k}\dfrac{(4k+3)c_m(2k,0)}{(8k+2)(k-m+1)}
	[(\frac{u}{t})^{2m}-(\frac{u}{t})^{2k+2}]\right)du
	\nonumber\\
	-\sum_{i=1}^{k+1}t^{2i-1}\int_{0}^{t}\sum_{m=1}^{k+i}\left([\frac{(4k+3)c_m(2k,2i)}{(8k+2)(k+i+1-m)}+\dfrac{s_m(2k+1,2(i-1)}{4(k+i-m+1)}]\right.\nonumber\\\left.\times[(\frac{u}{t})^{2m}-(\frac{u}{t})^{2k+2i+2}]+\dfrac{c_m(2k,2i)}{4k+1}(\frac{u}{t})^{2m}\right)(\Delta^{i}(Mf(p,q,u)))du.\,\,\,\,\,\,\,\label{5.14}
\end{eqnarray}

\noindent It follows from (\ref{5.14}) that (\ref{5.9}) is true for every $k$. Theorem 4 is proved.

\section{Algorithms to calculate the coefficients $c_m$ and $s_m$}\label{sec6}

\noindent Here we present recurrent algorithms to calculate the coefficients 
$c_m(2k+1,2i)$ for integers $k\geq1$, $0\leq i\leq k$, $1\leq m\leq k+i$ and the coefficients
$s_m(2k,2(i-1))$ for integers $k\geq1$, $0\leq i\leq k$, $1\leq m\leq k+i-1$, otherwise they are $0$. 

\noindent Comparing (\ref{5.3}) and (\ref{5.7}) we obtain: for $i=0$
\begin{equation}\label{6.1}\begin{cases}
		c_m(2k+1,0)=-(4k+3)\left(\frac{(4k+1)c_m(2k-1,0)}{(8k-2)(k-m)}-\frac{c_m(2k-1,0)}{4k-1}\right),\,\, 1\leq m\leq k-1\\
		c_{k}(2k+1,0)=-(4k+3)\left(\dfrac{4k+1}{3}- \sum_{m=1}^{k-1}
		\frac{(4k+1)c_m(2k-1,0)}{(8k-2)(k-m)}\right);\end{cases}
\end{equation}
for $1\leq i\leq k$
\begin{equation}\label{6.2}\begin{cases}
		c_m(2k+1,2i)=(4k+3)(\frac{c_m(2k-1,2i)}{4k-1}+\frac{-s_m(2k,2(i-1))}{4(k+i-m)}-\frac{(4k+1)c_m(2k-1,2i)}{(8k-2)(k+i-m)})\\\qquad\qquad\qquad\qquad\qquad\qquad\qquad\,\,\,\,\,\,\,\,\,\,\,\,\,\,\,\, 1\leq m< k+i\\
		c_{k+i}(2k+1,2i)=(4k+3)\sum_{m=1}^{k+i-1}\left(\frac{(4k+1)c_m(2k-1,2i)}{(8k-2)(k+i-m)}+\frac{s_m(2k,2(i-1)}{4(k+i-m)}\right).\end{cases}
\end{equation}

\noindent Comparing (\ref{5.4}) and (\ref{5.6})  we obtain: for $i=1$

\begin{equation}\label{6.3}\begin{cases}
		s_m(2k,0)=\frac{4k+1}{2k(2k+1)}\left(
		\frac{(2k-1)(2k-2)s_m(2(k-1),0)}{4k-3}+\frac{c_m(2k-1,0)}{k-m}-\right.\\\left.\frac{(2k-1)(k-1)(4k-1)s_m(2(k-1),0)}{(4k-3)(k-m)}\right)\,\, \texttt{for}\,\,\, 1\leq m\leq k-1\\
		s_{k}(2k,0)=\frac{4k+1}{2k(2k+1)}\left(\frac{8k-2}{3}-\sum_{m=1}^{k-1}\left(\frac{c_m(2k-1,0)}{k-m}-\right.\right.\\\left.\left.
		\frac{(2k-1)(k-1)(4k-1)s_m(2(k-1),0)}{(4k-3)(k-m)}\right)\right).
	\end{cases}
\end{equation}

for $1<i\leq k$, we get

\begin{equation}\label{6.4}\begin{cases}
		s_m(2k,2(i-1))=\left(\frac{(2k-1)(2k-2)s_m(2(k-1),2(i-1))}{4k-3}+\frac{c_m(2k-1,2(i-1))}{k+i-m-1}-\right.\\\left.\frac{(2k-1)(k-1)(4k-1)s_m(2(k-1),2(i-1))}{(4k-3)(k+i-m-1)}\right)\frac{4k+1}{2k(2k+1)}\,\, \texttt{for}\,\,\, 1\leq m\leq k+i-2\\
		s_{k+i-1}(2k,2(i-1))=-\frac{4k+1}{2k(2k+1)}\sum_{m=1}^{k+i-2}\left(\frac{c_m(2k-1,2(i-1))}{k+i-m-1}-\right.\\\left.
		\frac{(2k-1)(k-1)(4k-1)s_m(2(k-1),2(i-1))}{(4k-3)(k+i-m-1)}\right).
	\end{cases}
\end{equation}

\noindent Here we present recurrent algorithms to calculate the coefficients 
$c_m(2k,2i)$ for integers $k\geq1$, $0\leq i\leq k$, $1\leq m\leq k+i$ and the coefficients
$s_m(2k+1,2(i-1))$ for integers $k\geq1$, $1\leq i\leq k+1$, $1\leq m\leq k+i
$.

\noindent Comparing (\ref{5.14}) and (\ref{5.8}) we obtain: for $i=0$ 
\begin{equation}\label{6.5}\begin{cases}
		c_m(2k+2,0)=(4k+5)\left(\frac{c_m(2k,0)}{4k+1}-\frac{(4k+3)c_m(2k,0)}{(8k+2)(k-m+1)}\right), 1\leq m\leq k\\
		c_{k+1}(2k+2,0)=(4k+5)\left(\sum_{m=1}^{k}\frac{(4k+3)c_m(2k,0)}{(8k+2)(k-m+1)}-(4k+3)\right);\end{cases}
\end{equation}
and for $1\leq i\leq k+1$
\begin{equation}\label{6.6}\begin{cases}
		c_m(2(k+1),2i)=(4k+5)(\frac{-(4k+3)c_m(2k,2i)}{(8k+2)(k+i-m+1)}-\frac{s_m(2k+1,2(i-1))}{4(k+i-m+1)}+\frac{c_m(2k,2i)}{4k+1}), \\
		\,\,\,\,\,\,\,\,\,\,\,\,\,\,\,\text{for}\,\,\,\,\,1\leq m\leq k+i\\
		c_{k+i+1}(2(k+1),2i)=	(4k+5)\sum_{m=1}^{k+i}\left(\frac{(4k+3)c_m(2k,2i)}{(8k+2)(k+i-m+1)}+\frac{s_m(2k+1,2(i-1))}{4(k+i-m+1)}\right).
	\end{cases}
\end{equation}
\noindent Comparing (\ref{5.13}) and (\ref{5.9}) we obtain: for $i=1$

\begin{equation}\label{6.7}\begin{cases}
		s_m(2k+1,0)=\frac{4k+3}{(2k+2)(2k+1)}\left(\frac{2k(2k-1)s_m(2k-1,0)}{4k-1}+\frac{c_m(2k,0)}{k+1-m}\right.\\\left.-\frac{k(2k-1)(4k+1)s_m(2k-1,0)}{(4k-1)(k+1-m)}\right)\,\, \texttt{for}\,\,\, 1\leq m\leq k\\
		s_{k+1}(2k+1,0)=\frac{4k+3}{(2k+2)(2k+1)}\left(8k+2-\sum_{m=1}^{k}(\frac{c_m(2k,0)}{k+1-m}-\right.\\\left.\frac{k(2k-1)(4k+1)s_m(2k-1,0)}{(4k-1)(k+1-m)})\right).
	\end{cases}
\end{equation}

for $1<i\leq k+1$, we get

\begin{equation}\label{6.8}\begin{cases}
		s_m(2k+1,2(i-1))=\frac{4k+3}{(2k+2)(2k+1)}\left(\frac{2k(2k-1)s_m(2k-1,2(i-1))}{4k-1}+\right.\\\left.\frac{c_m(2k,2(i-1))}{k+i-m}-\frac{k(2k-1)(4k+1)s_m(2k-1,2(i-1))}{(4k-1)(k+i-m)}\right), 1\leq m< k+i\\
		s_{k+i}(2k+1,2(i-1))=-\frac{4k+3}{(2k+2)(2k+1)}\left(\sum_{m=1}^{k+i-1}(\frac{c_m(2k,2(i-1))}{k+i-m}-\right.\\\left.\frac{k(2k-1)(4k+1)s_m(2k-1,2(i-1))}{(4k-1)(k+i-m)})\right).
	\end{cases}
\end{equation}

\noindent From (\ref{5.5}) - (\ref{5.6}) we get
\begin{equation}\label{6.9}
	c_1(2,0)=-15; \,\,c_2(2,0)=0; \,\,s_1(1,0)=3.
\end{equation}
Finally, using (\ref{6.9}) and  the algorithms  one can calculate all coefficients

\noindent $c_m(2(k+1),2i)$ $s_m(2k,2i)$ for integers $k\geq1$, $0\leq i\leq k$, $\,1\leq m\leq k+i$ and 

\noindent $c_m(2k+1,2i)$ $s_m(2k+1,2i)$ for integers $k\geq1$,$\,0\leq i\leq k$, $\,1\leq m\leq k+i$.

\section{Proof of (\ref{1.6})}\label{sec7}

\noindent Using Theorem 3 and Theorem 4 for the partial sum of (\ref{3.3}), we get

\begin{eqnarray}
	\sum_{k=0}^{n} (a_{0(2k)}(x,y,\mid z\mid)+
	(\sgn{\,z})a_{0(2k+1)}(x,y,\mid z\mid))=(n+1)\times\,\,\,\,\,\label{7.1}\\
	\left((2n+1) Mf(x,y,\mid{z}\mid)+(\sgn{\,z})\dfrac{(2n+3)}{3}a_{01}(x,y,\mid z\mid)\right)+\sum_{i=0}^{n}\int_{0}^{\mid z\mid}\mid z\mid^{2i-1}\nonumber\\\times
	\left(C_{2n,2i}
	(\dfrac{u}{\mid z\mid})\Delta^{i}(Mf(x,y,u))+ (\sgn{z})C_{2n+1,2i}
	(\dfrac{u}{\mid z\mid})\Delta^{i}(a_{01}(x,y,u))\right)du\nonumber\end{eqnarray}
where $C_{2n,2i},$ and $C_{2n+1,2i},$ are the following polynomials $0\leq i\leq n$:
\begin{eqnarray}\label{7.2}
	C_{2n,2i}(t)= \sum_{k=i}^{n}\sum_{m=1}^{k+i}c_m(2k,2i)\,t^{2m},\\
	C_{2n+1,2i}(t)= \sum_{k=i}^{n}\sum_{m=1}^{k+i}c_m(2k+1,2i)\,t^{2m+1}.\nonumber\end{eqnarray}
Here, we assume that $C_{0,0}\equiv0$ and $C_{1,0}\equiv0$.
\noindent Substituting (\ref{7.1}) into  (\ref{3.3}), we obtain (\ref{1.6}). Theorem 1.2 is proved.

\section{Computational Implementation}\label{sec8}

\noindent Now we are going to implement the  iterative reconstruction algorithm obtained from (\ref{1.6}).

\noindent	\textbf{Example}. Consider the following odd function (bad case from SRT point of view) with respect to the $z=0$ plane 

\begin{equation}\label{8.1}
	f(x,y,z)=			(x^2+y^2) z^3.
\end{equation}

\noindent Fig. 1 shows the plot (using Mathematica) of the restriction of $f$ (onto the plane $y=3$) on $[-3,3]\times[-2,2]$.

\begin{figure}
	\center
	\includegraphics{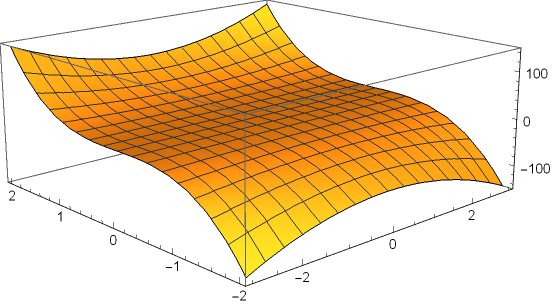}
	\caption{the restriction of $f$ (onto the plane $y=3$) on $[-3,3]\times[-2,2]$.}
\end{figure}

\noindent Now we approximate $f$ for $n=2$ using the reconstruction algorithm obtained from (\ref{1.6}).
For $Mf(x,u)$ and $a_{01}(x,u)$ of the restriction of $f$ over the sphere with the center
$(x,y,0)$ and radius $u\geq 0$, we get 	\begin{equation}\label{8.2}
	Mf(x,y,u)\equiv0  \,\,\,\,\emph{and}\,\,\,\,  a_{01}(x,y,u)=\frac{3}{5}(x^2+y^2)u^3 + \frac{6}{35}u^5.
\end{equation}
Using (\ref{6.1})-(\ref{6.4}) and (\ref{7.2}), for $n=2$, we get
\begin{equation}\label{8.3}
	C_{5,0}(t)=\frac{105}{2}t^3-\frac{7\cdot 33}{2}t^5\,\,\,\,\emph{and}\,\,\,\, C_{5,2}(t)=\frac{75\cdot 5}{16}t^3+\frac{219}{16}t^5-\frac{11\cdot 33}{16}t^7.	\end{equation}

\noindent Substituting (\ref{8.2}) and (\ref{8.3}) into (\ref{1.6}) we obtain the approximation of $f$ for $n=2$. Fig. 2 shows the plot of the approximation of $f$ on $[-3,3]\times[-2,2]$ (on plane $y=3$).

\begin{figure}
	\center
	\includegraphics{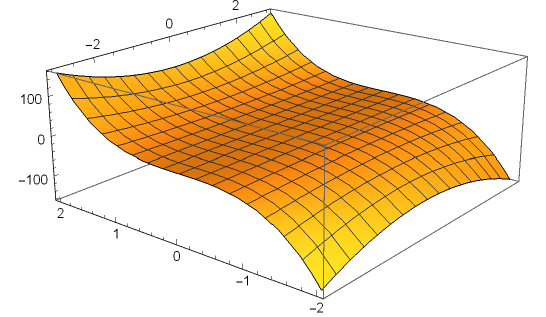}
	\caption{the approximation of $f$ obtained by (\ref{1.6}) for $n=2$ on $[-3,3]\times[-2,2]$ ($y=3$).}
\end{figure}

\end{document}